\documentclass{amsart}

\usepackage{amsmath,amsthm,amssymb,amscd,enumerate,mathtools,enumitem}
\usepackage{amsfonts}
\usepackage{rotating}
\usepackage{euscript}
\usepackage{pst-node}
\usepackage{epsfig,verbatim}

\def\be{\begin{equation}}
\def\ee{\end{equation}}

\def\C{{\mathbb C}} 
\def\f{\EuScript}
\def\N{{\mathbb N}}

\def\ord{{\rm ord\,}}

\def\phi{{\varphi}}
\def\v{{\varepsilon}}

\def\bp{\begin{proposition}}
\def\ep{\end{proposition}}

\def\bt{\begin{theorem}}
\def\et{\end{theorem}}
\def\br{\begin{remark}}
\def\er{\end{remark}}
\def\be{\begin{equation}}
\def\bee{\begin{equation*}}
\def\l{\label}

\def\ee{\end{equation}}
\def\eee{\end{equation*}}
\def\bl{\begin{lemma}}
\def\el{\end{lemma}}
\def\bc{\begin{corollary}}
\def\ec{\end{corollary}}
\def\pr{\noindent{\it Proof. }}

\def\bd{\begin{definition}}
\def\ed{\end{definition}}

\newtheorem{theorem}{Theorem}[section]
\newtheorem{lemma}[theorem]{Lemma}
\newtheorem{definition}[theorem]{Definition}
\newtheorem{corollary}[theorem]{Corollary}
\newtheorem{proposition}[theorem]{Proposition}
\newtheorem{problem}[theorem]{Problem}

\theoremstyle{definition}

\theoremstyle{definition}
\newtheorem{remark}[theorem]{Remark}
\def\bpr{\begin{problem}}
\def\epr{\end{problem}}

\mathtoolsset{showonlyrefs}

\begin{document}
\title{Right amenability in semigroups of formal power series}
\author{Fedor Pakovich}
\thanks{
This research was supported by ISF Grant  No. 1092/22}
\address{Department of Mathematics, Ben Gurion University of the Negev, Israel}
\email{
pakovich@math.bgu.ac.il}

\begin{abstract}
Let $k$ be an algebraically closed field of characteristic zero, and $k[[z]]$ the  ring of formal power series over $k$. We provide several characterizations of right amenable    finitely generated subsemigroups of  $z^2k[[z]]$ with the semigroup operation $\circ $ being composition. In particular,  
we show that a subsemigroup $S=\langle Q_1,Q_2,\dots, Q_k\rangle$ of  $z^2k[[z]]$ is right amenable  if and only if there exists an invertible element $\beta$ of $zk[[z]]$
such that $\beta^{-1}\circ Q_i \circ \beta =\omega_i z^{d_i},$ $1\leq i \leq k,$ for some integers $d_i$, $1\leq i \leq k,$ and roots of unity $\omega_i,$ $1\leq i \leq k.$

\end{abstract}

\maketitle

\section{Introduction} 

Let $R$ be a commutative  ring with identity, and $R[[z]]$ the ring of formal power series over $R$. For an element $A(z)=\sum_{n\geq 0}c_nz^n$ of $R[[z]]$,  its {\it order} is defined   
by the formula $\ord A=\min\{n\geq 0\,\vert\, c_n\neq 0\}.$  
If $A$ and $B$ are elements of $R[[z]]$ with \linebreak $\ord B\geq 1$, then the operation $A\circ B$ of  composition of $A$ and $B$ is well defined and provides $zR[[z]]$ with the structure of a  semigroup. This semigroup 
contains a  group $\mathcal J(R)$ consisting  of all formal power series of the form $z+\sum_{n\geq 2}c_nz^n$.  The group $\mathcal J(R)$  has been  extensively studied (see survey \cite{bab1} and references therein).   
In particular, it was established in \cite{bab2} that $\mathcal J(R)$ is  amenable as a {\it topological} group.   In this note, we study the right amenability of  subsemigroups of $z^2k[[z]]$, where  $k$ denotes  an algebraically closed field of characteristic zero. However, in distinction with \cite{bab2}, all studied semigroups are considered as {\it discrete}. 
This setting is different,  and requires another approach. In a sense, 
the  results of this note can be seen as analogues of the results of the recent papers \cite{peter}, \cite{peter2}, \cite{amen} about right amenable semigroups of polynomials and rational functions over $\C$. Nevertheless,  our  methods are  different.

Let us recall that a semigroup $S$ is called {\it right amenable} if it admits a finitely additive probability measure $\mu$ defined on all the subsets of $S$  
such that for all $a\in S$ and $T\subseteq S$ 
the equality \be \l{perd0} \mu(Ta^{-1}) = \mu(T)\ee holds, where   
 the set $Ta^{-1}$ is defined by the formula  $$ Ta^{-1}=\{s \in S\, | \,sa \in T\}.$$ 
 A semigroup $S$ is called {\it right reversible}  if for all $a,b\in S$ the left ideals 
$Sa$ and  $Sb$ have a non-empty intersection, that is, if for all $a,b\in S$  there exist $x,y\in S$ such that 
$ xa=yb.$ 
It is well-known and follows easily from the definition  that any right amenable semigroup is right reversible.

We denote by  $\f Z^U$ the subsemigroup of $z^2k[[z]]$ consisting of all monomials  of 
the form $\omega z^n,$ $n\geq 2,$ where $\omega$ is a root of unity.
We say that two subsemigroups $S_1$ and $S_2$ of $z^2k[[z]]$ are {\it conjugate} if 
there exists   $\beta\in k[[z]]$ of order one 
 such that $$\beta^{-1}\circ S_1\circ \beta=S_2.$$ 
In this notation, our main result is following.

\bt \l{main0} Let $k$ be an algebraically closed field of characteristic zero, and $Q_1,Q_2,\dots, Q_k$  elements of $z^2k[[z]]$. Then for the semigroup  $S=\langle Q_1,Q_2,\dots, Q_k\rangle$  generated by  $Q_1,Q_2,\dots, Q_k$  the following conditions are equivalent:

\begin{enumerate} [label=\normalfont \arabic*)]

\item The semigroup $S$ is right amenable. 

\item The semigroup $S$ is right reversible. 

\item The semigroup $S$ contains  no free subsemigroup of rank two.  

\item  The intersection of principal left ideals    $SQ_1\cap SQ_2\cap \dots \cap SQ_k$ is non-empty.

\item The semigroup $S$ is conjugate to a subsemigroup of $\f Z^U$.

\end{enumerate}

\et 

The rest of this note is organized as follows. In Section 2.1,  we collect some auxiliary results that are used in the paper. In Section 2.2,  we show that  every subsemigroup $S$ of $\f Z^U$, not necessarily finitely generated,  is right amenable and contains  no free subsemigroup of rank two. Then, using the result of the  paper \cite{sha}, which is essentially equivalent to the equivalence $4)\Leftrightarrow 5)$ in Theorem \ref{main0}, we prove Theorem \ref{main0}.  Finally, in Section 2.3, we provide a class of examples showing that Theorem \ref{main0} is not true for  {\it infinitely} generated subsemigroups of $z^2k[[z]]$.

\section{Proof of Theorem \ref{main0}}

\subsection{Auxiliary results} Let us recall that a semigroup $S$ is called {\it left cancellative} if the equality $ab=ac,$ where $a,b,c\in S,$ implies that $b=c$. Right cancellative semigroups are defined similarly. A  semigroup $S$ is called {\it cancellative} if it is left and right cancellative. 

The following result relates left reversibility with the presence of  free subsemigroup of rank two (see \cite{frey}, Theorem 8.4, or \cite{gray}, Corollary 4.2).

\bl \l{xrun}  Let $S$ be a right cancellative semigroup that contains no free subsemigroup of rank two. Then $S$ is right reversible. \qed 
\el

A subsemigroup of a right amenable
semigroup  is not necessarily right amenable. However, the following result holds (see \cite{frey}, Theorem 8.5, or \cite{don}, Theorem 4). 

\bt \l{fr} Let $S$ be a cancellative semigroup such that $S$ contains no free
subsemigroup on two generators. If $S$ is right amenable, then every subsemigroup
of $S$ is right amenable. \qed 
\et 
   
For a semigroup $U$, we denote by  ${\rm End}(U)$ the set of endomorphisms of $U$. 
Suppose that $U$ and $T$ are semigroups with a homomorphism
$\rho: T \rightarrow {\rm End}(U).$  Denoting for $a\in T$ the endomorphism $\rho(a)$ of $U$  by $\rho_a,$  we define the semidirect product of $U$ and $T$
as the semigroup $F=U\underset{\rho}{\times} T$ of ordered pairs $(u, a)$, where 
$u\in U$ and $a\in T$, with the operation 
$$(u, a)(v, b) = (u\rho_a(v), ab).$$ 

An example of a semidirect product is provided by the  subsemigroup $\f Z$  of $z^2k[[z]]$  consisting of all   monomials  $az^n,$ $n\geq 2,$ where $a\in k^*$. 
Indeed, for every integer $n\geq 2$ the map $\phi_{n}:\ a\rightarrow a^{n}$ is
an endomorphism of $k^*$. Moreover, the corresponding map 
$n\rightarrow \phi_{n}$ induces a  homomorphism 
\be \l{hom} \rho: \N\setminus \{1\} \rightarrow {\rm End}(k^*),\ee where $\N\setminus \{1\}$ stands for the corresponding subsemigroup of the multiplicative semigroup of natural numbers, identified with the subsemigroup  $z^n,$ $ n\geq 2$, of $z^2k[[z]].$     Thus, the 
semidirect product $k^*\underset{\rho}{\times} \N\setminus \{1\}$ is well defined and can be identified with $\f Z$. 
Similarly, $\f Z^U$ can be identified with the 
semidirect product $\mu_{\infty} \underset{\rho}{\times} \N\setminus \{1\}$, where  $\mu_{\infty}$ is the group of all roots of unity. 
More generally, we can define the 
semidirect product $U\underset{\rho}{\times} N$, where $U$ is any subsemigroup of $k^*$, $N$ is any  subsemigroup $\N\setminus\{1\}$, and $\rho$ is the restriction of the homomorphism \eqref{hom} on $N.$ Below, we will omit  $\rho$ in the notation of such semigroups.

The following result was proved in \cite{kla}. 

\bt \l{kla} If $U$ and $T$ are right amenable semigroups
with a homomorphism $\rho: T \rightarrow {\rm End}(U)$, then $F=U\underset{\rho}{\times} T$ is right
amenable. 
\qed 
\et 
 
Let us recall that a {\it congruence} on a semigroup $S$ is an equivalence relation on $S$ compatible with the structure of semigroup, that is, an equivalence relation such that $x\sim y$ and $x'\sim y'$ implies that $xx'\sim yy'.$ 
If $S$ is  a semigroup and $\sim$ is a congruence on $S$, then one can define the {\it quotient
semigroup} $S/\sim$, whose elements are the equivalence  classes of $\sim$, and for $a, b \in S$ the operation on the corresponding classes  is defined by $[a] * [b] = [ab].$ 

Congruences correspond to homomorphic images of $S$ in the following sense.  If $\phi:\, S\rightarrow T$ is a homomorphism of semigroups, then the equivalence relation $\sim_{\phi} $, defined by  $x\sim_{\phi} y$ if and only if $\phi(x)=\phi(y)$, is a congruence on $S$ and the isomorphism 
\be \l{isom} 
 \phi(T)\cong S/\sim_{\phi}
\ee
holds.

Let $S$ be a right reversible semigroup, and let $\sim$ be the relation on $S$, defined by $x\sim y$ if and only if there exists $s\in S$ such that \be \l{listik} s\circ x=s\circ y.\ee In this notation, the following criterion for the right amenability holds (see  \cite{pater}, Proposition 1.24 and Proposition 1.25).

\bt \l{lastik} Let $S$ be a right reversible semigroup. 
Then the relation $\sim$ is a congruence on $S$, and the semigroup $S/\sim$ is left cancellative. Moreover, $S$  is right amenable if and only if $S/\sim$ is right amenable. \qed 
\et

Let $A\in z^2k[[z]]$ be 
 a formal power series of order $n$.  We recall that a {\it B\"ottcher function} associated with $A$ is a formal series  $\beta_A$ of order one  
such that  the equality 
\be \l{a} A\circ \beta_A=\beta_A\circ z^n\ee holds. 
 It is known that such a function exists and  is defined in a unique way 
up  to the change $\beta_A(z)\rightarrow \beta_A(\v z),$ where $\v^{n-1}=1$ 
(see \cite{kau}, Hilffsatz 4). Among other things, 
it follows from the existence of B\"ottcher functions that the semigroup $z^2k[[z]]$ is right cancellative. Indeed, if 
\be \l{b0} 
A_1\circ X=A_2\circ X,
\ee 
then conjugating  \eqref{b0} by $\beta_X$ we obtain the equality 
$$\beta_X^{-1} \circ A_1 \circ \beta_X\circ z^n=\beta_X^{-1} \circ A_2 \circ \beta_X\circ z^n, 
$$
which implies that 
$$\beta_X^{-1} \circ A_1 \circ \beta_X=\beta_X^{-1} \circ A_2 \circ \beta_X$$ and $A_1=A_2.$

Finally, we need the following result.  

\bt \l{main1}  Let $k$ be an algebraically closed field of characteristic zero, \linebreak $Q_1,Q_2,\dots, Q_k$   elements of $z^2k [[z]]$, and  $S=\langle Q_1,Q_2,\dots, Q_k\rangle$  the semigroup  generated by $Q_1, Q_2,\dots, Q_k$. Assume that $Q_1$ is contained in $\f Z^U$. Then  
\be \l{xera} SQ_1\cap SQ_2\cap \dots \cap SQ_k\neq\emptyset\ee
if and only if every $Q_i,$  $2\leq i \leq k$, is contained in $\f Z^U$. \qed
\et
 
In case $k=\C$, Theorem \ref{main1} was proved in \cite{sha} (Theorem 2.3),  
and the proof carries over verbatim to the case of an arbitrary algebraically closed field $k$  of characteristic zero.

\subsection{Right amenability of subsemigroups  $\f Z^U$.}
Let us denote by $\mu_n$ the   group of  $n$th roots of unity. 

\bt \l{m1} Every subsemigroup $S$ of $\f Z^U$ is right reversible and contains  no free subsemigroup of rank two.  
\et 
\pr Let us show that $S$ contains  no free subsemigroup of rank two.  Let 
$$F_1=\v_1z^{d_1}, \ \ \ \  F_2=\v_2z^{d_2}, \ \ \ \  \v_1,\v_2\in \mu_{\infty},\ \ \ \ d_1,d_2\geq 2,$$ be elements of $S$, and $n\geq 1$  an integer such that $\v_1,\v_2\in \mu_n.$ Assume first that $d_1=d_2.$ 
Then for every $j\geq 1$ the equality  
$$F_1^{\circ j}=\omega_jF_2^{\circ j}$$  holds for some $\omega_j\in \mu_n$, implying by the pigeonhole principle that there exist $j_1\neq j_2$  such that  
$$F_1^{\circ j_1}=\v F_2^{\circ j_1}, \ \ \   F_1^{\circ j_2}=\v F_2^{\circ j_2}$$ for the same $\v\in \mu_n.$ Assuming that $j_2>j_1$, this yields that   
\be \l{u1} F_1^{\circ j_2}=F_1^{\circ j_1}\circ F_2^{\circ (j_2-j_1)},\ee and hence the semigroup $<F_1,F_2>$ 
generated by $F_1$ and $F_2$ is not free.

In case $d_1\neq d_2,$ let us consider the elements  $$F_1'=F_1\circ F_2, \ \ \ \ F_2'=F_2\circ F_1$$ of $S$. If $F_1'=F_2'$, then  the semigroup $<F_1,F_2>$ obviously is not free. 
On the other hand, if $F_1'\neq F_2',$ then, since $F_1'$ and $F_2'$ have the same order, 
we can apply the above reasoning to $F_1'$ and $F_2'$,  and find $j_2>j_1$  such that   
\be \l{u2} (F_1\circ F_2)^{\circ j_2}=(F_1\circ F_2)^{\circ j_1}\circ (F_2\circ F_1)^{\circ (j_2-j_1)}.\ee

To prove that $S$ is right reversible it is enough to observe that equalities \eqref{u1} and \eqref{u2} provide solutions 
of  the equation \be \l{th} X\circ F_1=Y\circ F_2 \ee  in $X,Y\in S.$  \qed 

\vskip 0.2cm
For $\v\in \mu_{\infty}$, we denote by $\vert \v \vert$ the order of $\v$ in the semigroup 
$\mu_{\infty}.$ With every subsemigroup $S$ of  $\f Z^U$ we associate several objects. First, we define $U(S)$ as the subsemigroup of $\mu_{\infty}$ generated by all roots of unity $\v$ such that $\v z^d\in S$ for some  $d\geq 2.$ Notice that since $\v^{-1}=\v^{\vert \v\vert -1}$ belongs to  $U(S)$ whenever $\v$ belongs to $ U(S)$, 
the semigroup $U(S)$ is a group. Second, we define $N(S)$ as   the subsemigroup of the multiplicative semigroup $\N\setminus\{1\}$ consisting of all $d\geq 2$ such that $\v z^d \in S$  for some $\v\in \mu_{\infty}.$ Notice that  by construction the semigroup $S$ is a subsemigroup of the semigroup  $U(S)\times N(S)$.

Further, we associate with $S$ two subsets $P_1(S)$ and $P_2(S)$ of the set of prime numbers as follows. For an integer $l\geq 2$, we define $\f P(l)$ as the set of prime divisors of $l.$  
For an element $Q= \v z^d$ of $S$,    where     $\v\in \mu_{\infty}$  and $d\geq 2$, we set $$p_1(Q)=\f P(\vert \v \vert ), \ \ \ \ \ p_2(Q)=\f P(d)$$
(in case $\vert \v \vert =1$, we set $p_1(Q)=\emptyset$). Finally, we set
$$P_1(S)=\bigcup_{Q\in S}p_1(Q), \ \ \ \  \  P_2(S)=\bigcup_{Q\in S}p_2(Q).$$ 
Notice that  by construction 
\be \l{by} P_1(S)=P_1\left(U(S)\times N(S)\right), \ \ \  \ P_2(S)=P_2\left(U(S)\times N(S)\right).\ee

\bt \l{m2} Every subsemigroup $S$ of $\f Z^U$ is right amenable. 
\et 
\pr 
We start by observing that  any semigroup of the form $U{\times} N$, where $U$ is a subsemigroup of $k^*$ and $N$ is a subsemigroup $\N\setminus\{1\}$ is  right amenable. 
Indeed, it is well known that every commutative semigroup is right amenable. Therefore, 
  $U$ and $N$ are right amenable, implying by Theorem \ref{kla} that $U{\times} N$ is also right amenable.

Let us show first 
that the statement of the theorem holds if \be \l{usl} P_1(S)\cap P_2(S)=\emptyset.\ee Since $S$ is  a subsemigroup of the amenable semigroup $U(S)\times N(S)$ and the last group is right cancellative and  contains  no free subsemigroup of rank two by Theorem \ref{m1}, it follows from Theorem \ref{fr} that to prove the amenability of $S$ it is enough to show that the semigroup $U(S)\times N(S)$ is left cancellative. 
Let us assume that 
\be \l{bur} A\circ F_1=A\circ F_2\ee for some $A,F_1,F_2\in U(S)\times N(S)$. Clearly, \eqref{bur} implies that $\ord F_1=\ord F_2$. Thus, 
$$A=\v z^m, \ \ \ F_1=\v_1 z^k, \ \ \  F_2=\v_2 z^k$$ for some $m,k\in N(S)$ and $\v,$ $\v_1,$ $\v_2\in U(S),$ and 
 \eqref{bur} implies the equality 
\be \l{ip} \left(\frac{\v_1}{\v_2}\right)^m=1.\ee 
Set  $n=\vert \v_1/\v_2\vert.$ Since \eqref{by} and  \eqref{usl} yield that $\gcd(n,m)=1,$ 
equality \eqref{ip} is possible only if $n=1$. Therefore,   
$\v_1=\v_2$ and $F_1=F_2.$ 

Let us prove now the theorem  in the general case. 
Since for any integers $n,m_1,m_2\geq 1$ such that $n=m_1m_2$ and  $\gcd(m_1,m_2)=1$ the equality 
$$\mu_n=\mu_{m_1}\times \mu_{m_2}$$ holds, any element $Q$ of $S$ has a unique representation 
in the form 
\be \l{def} Q=\v_1\v_2z^{d},\ \ \ \ \v_1,\v_2\in \mu_{\infty},\ \ \ \ d\geq 2,\ee
where $$\f P(\vert \v_1\vert)\cap    P_2(S)=\emptyset, \ \ \ \ \f P(\vert \v_2\vert)\subseteq   P_2(S).$$
Let us define a map 
$$\phi:\, S\rightarrow \f Z^U,$$ setting for $Q$, defined by \eqref{def},  
$$\phi(Q)=\v_1z^{d}.$$ It is easy to see that $\phi$ is a semigroup homomorphism. 
Indeed, for  
$$\widehat Q=\widehat\v_1\widehat\v_2z^{\widehat d},\ \ \ \ \widehat\v_1,\widehat\v_2\in \mu_{\infty},\ \ \ \ \widehat d\geq 2,$$ 
 we have:
$$Q\circ \widehat Q=\v_1\v_2z^{d}\circ \widehat\v_1\widehat\v_2z^{\widehat d}=\v_1{\widehat\v_1}^{\ d}\v_2{\widehat\v_2}^{\ d}z^{d\widehat d}, 
$$ where obviously 
$$\f P(\vert \v_1{\widehat\v_1}^{\ d}\vert )\cap    P_2(S)=\emptyset, \ \ \ \  \f P(\vert \v_2{\widehat\v_2}^{\ d}\vert)\subseteq   P_2(S).$$ 
Thus, 
$$\phi(Q\circ \widehat Q)=\v_1{\widehat\v_1}^{\ d}z^{d\widehat d}=\phi(Q)\circ \phi(\widehat Q).$$ 
By construction, the image $\phi(S)$ of $S$ under the homomorphism $\phi$ satisfies the condition 
$$P_1(\phi(S))\cap P_2(\phi(S))=\emptyset.$$ Therefore, by what is proved above, the semigroup  $\phi(S)$ is cancellative and right amenable. 
It follows now from Theorem \ref{lastik} taking into account the isomorphism \eqref{isom} that to prove the theorem 
it is enough to show that the congruence defined by the homomorphism $\phi$ coincides with the congruence \eqref{listik}.

Let us assume that 
\be \l{hi} \phi(Q_1)=\phi(Q_2).\ee 
Then $$Q_1=\v\v_1z^{d}, \ \ \ \ Q_2=\v\v_2z^{d}, \ \ \ \ \v,\v_1,\v_2\in \mu_{\infty},\ \ \ \ d\geq 2,$$ where  
 $$\f P(\vert \v\vert)\cap    P_2(S)=\emptyset, \ \ \ \ \f P(\vert \v_1\vert)\subseteq   P_2(S),\ \ \ \ \f P(\vert \v_2\vert)\subseteq   P_2(S).$$
Let 
\be \l{ass} \vert \v_1\vert=p_1^{a_1}p_2^{a_2}\dots p_r^{a_r}, \ \ \ \ \ \vert \v_2\vert=q_1^{b_1}q_2^{b_2}\dots q_l^{b_l}\ee be the canonical decompositions of
$\vert \v_1\vert$ and $ \vert \v_2\vert$  into products of primes. Then for each $i,$ $1\leq i \leq r$, there exists $K_i\in S$ of order divisible by $p_i$, and for each $j,$ $1\leq j \leq l$, there exists $L_j\in S$ of order divisible by $q_j$. Obviously, this implies that the equality 
\be \l{ih} A\circ Q_1=A\circ Q_2\ee holds for the element 
$$A=K_1^{\circ a_1}\dots K_r^{\circ a_r}L_1^{\circ b_1}\dots L_l^{\circ b_l}$$
 of $S$  (formulas  \eqref{ass} assume that $\vert \v_1\vert>1$, $ \vert \v_2\vert>1$, however, in case one of these numbers equals one  the proof  can be  modified in an obvious way).  

In the other direction, if \eqref{ih} holds, then 
$$\phi(A)\circ \phi(Q_1)=\phi(A)\circ \phi(Q_2),$$ implying that \eqref{hi} holds since 
 $\phi(S)$ is cancellative.  Thus, the congruence defined by the homomorphism $\phi$ coincides with the congruence \eqref{listik}.  \qed

\vskip 0.2cm 
\noindent {\it Proof of Theorem \ref{main0}.} We first prove the chain of implications 
$1)\Rightarrow 2)\Rightarrow 4)\Rightarrow 5)\Rightarrow 1).$
It is well known (see \cite{pater}, Proposition 1.23) that every right amenable semigroup is right reversible. Thus, $1)\Rightarrow 2).$ The implication $2)\Rightarrow 4)$ is proved by induction on $k.$ For $k=2$, the condition  $2)$  coincides with the right reversibility condition. On the other hand, if 
$$F=A_1\circ Q_1=A_2\circ Q_2=\dots =A_{k-1}\circ Q_{k-1}$$ for some $A_1,A_2,\dots A_{k-1}\in S$, then applying the condition of right reversibility 
to $F$ and $Q_k$, we can find $G,Y_1,Y_2\in S$  such that the equality 
$$G=Y_1\circ F=Y_2\circ Q_k$$ holds. Thus, 
$$G=(Y_1\circ A_1)\circ Q_1=(Y_1\circ A_2)\circ Q_2=\dots =(Y_1\circ A_{k-1})\circ Q_{k-1}=Y_2\circ Q_k.$$ 

The implication $4)\Rightarrow 5)$ follows from Theorem \ref{main1}.  Indeed, if $\beta_{Q_1}$
is a B\"ottcher function corresponding to $Q_1$, then 
$\beta_{Q_1}^{-1}\circ S\circ \beta_{Q_1}$ is a semigroup satisfying the conditions $Q_1\in \f Z^U$ and \eqref{xera}. Thus, by Theorem \ref{main1},  every $Q_i,$  $2\leq i \leq k$, is contained in $\f Z^U$, and hence  $S$ is conjugate to a subsemigroup of $\f Z^U$. Finally, the  implication $(5)\Rightarrow (1)$ follows from Theorem \ref{m2}. 

To finish the proof it is enough to prove the implications $5)\Rightarrow 3)$ and $3)\Rightarrow 2)$. The first implication follows from Theorem \ref{m1}. On the other hand, since 
$z^2k[[z]]$ is right cancellative, the second implication follows from Lemma \ref{xrun}. \qed

\subsection{Infinitely generated semigroups.}

Theorem \ref{main0} is not true for infinitely generated subsemigroups of $z^2k[[z]]$. 
A class of counterexamples is provided by Theorem \ref{last} below.

We recall that an element of a semigroup $S$ is called {\it indecomposable} if it belongs 
to $S\setminus SS,$ where $$SS = \{st : s, t \in S\}.$$

\bt \l{last} Let  $U$ be a subsemigroup of $k^*$ such that $U\not\subseteq \mu_{\infty}$ and $1\in U$. Then for any subsemigroup $N$  of $\N\setminus\{1\}$ the semigroup $S=U{\times} N$ satisfies the following conditions: 

\begin{enumerate} [label=\normalfont \arabic*)]

\item The semigroup $S$ is right amenable. 

\item The semigroup $S$ is not finitely generated. 

\item   The semigroup $S$ is not conjugate to a subsemigroup of $\f Z^U$.

\end{enumerate}

\et 
\pr The right amenability of $S$ was established in the proof of Theorem \ref{m2}. 
Further, since indecomposable elements of $S$ must belong to any generating
set of $S$, to prove that $S$ is not finitely generated it is enough to find an infinite subset of  indecomposable elements of $S$. An example of 
such a set is the set 
$a^kz^d,$ $ k\geq 1,$ where $a$ is an arbitrary element of $U$ that does not belong to $\mu_{\infty}$ and $d$ is an indecomposable element of the semigroup $N.$ 

 Finally, the proof of Theorem \ref{m1} shows that to prove that $S$ is not conjugate to a subsemigroup of $\f Z^U$ it is enough to find  $F_1,F_2\in S$ such that $\ord F_1=\ord F_2$ but equality \eqref{u1} does not hold for any choice of $j_1,j_2$. Since $1\in U,$  we can take  
$$F_1=z^d, \ \ \ \ F_2=az^d,$$ where $d$ is an arbitrary element of $N$, and $a$ is an element of $U$ that does not belong to $\mu_{\infty}.$ \qed


\begin{thebibliography}{10}

\bibitem{bab1} I. Babenko, {\it  Algebra, geometry and topology of the substitution group of formal power series}, Russian Math. Surveys 68 (2013), no. 1, 1-68 

\bibitem{bab2} I. Babenko, S. Bogatyi, {\it Amenability of the substitution group of formal power series}, Izv. Math. 75 (2011), no. 2, 239-252.


\bibitem {peter}  Cabrera C., Makienko P., {\it Amenability and measure of maximal entropy for semigroups of rational map},  Groups Geom. Dyn. 15 (2021), no. 4, 1139-1174.



\bibitem {peter2}  Cabrera C., Makienko P., {\it Amenability and measure of maximal entropy for semigroups of rational map: II}, arXiv:2109.11601.




\bibitem {don} J. Donnelly, {\it Subsemigroups of cancellative amenable semigroups}, Int. J. Contemp. Math. Sci. 7 (2012), no. 21-24, 1131-1137.

\bibitem {frey} Frey, A. H., {\it Studies on amenable semigroups}, 
Thesis (Ph.D.)–University of Washington. 1960. 




\bibitem{gray} R. D. Gray, M. Kambites, {\textit Amenability and geometry of semigroups},  Trans. Amer. Math. Soc. 369 (2017), no. 11, 8087-8103.




\bibitem {kau} H. Kautschitsch, \textit{\"Uber vertauschbare Potenzreihen,} Math. Nachr. 88 (1979), 207-217.



\bibitem{kla}  M. Klawe, \textit{
Semidirect product of semigroups in relation to amenability, cancellation properties, and strong Følner conditions}, 
Pacific J. Math. 73 (1977), no. 1, 91-106.







\bibitem{sha} F. Pakovich, {\it  Sharing a measure of maximal entropy in polynomial semigroups,} 
 Int. Math. Res. Not. IMRN 2022, no. 18, 13829-13840.

\bibitem{amen} F. Pakovich, {\it On amenable semigroups of rational functions,} 
Trans. Amer. Math. Soc. 375 (2022), no. 11, 7945-7979. 



\bibitem{pater} A. Paterson, {\it Amenability}, Mathematical Surveys and Monographs, 29. American Mathematical Society, Providence, RI, 1988. 




\end{thebibliography}
\end{document}